\documentclass[a4,11pt]{article}
\usepackage{amsfonts}
\usepackage{amsmath}
\usepackage{amssymb}
\usepackage{amsthm}
\theoremstyle{plain}
 \newtheorem{theorem}{Theorem}[section]
 
 \newtheorem*{thghjkl: eorem*}{Theorem}
 
 \newtheorem{lemma}[theorem]{Lemma}
 \newtheorem{corollary}[theorem]{Corollary}

\theoremstyle{definition}
 \newtheorem{definition}[theorem]{Definition}
\theoremstyle{remark}

\numberwithin{equation}{section}

\newcommand{\nexteq}{\displaybreak[0]\\ &=}
\newcommand{\nnexteq}{\nonumber\displaybreak[0]\\ &=}
\newcommand{\ts}{\theta^*}

\title{Coherent configurations and triply regular association schemes obtained from spherical designs}

\author{Sho Suda\\
{\small Division of Mathematics, Graduated School of Information Sciences, Tohoku University,} \\
{\small 6-3-09 Aramaki-Aza-Aoba, Aoba-ku, Sendai 980-8579, Japan}}
\date{\today}
\begin{document}
\maketitle

\begin{abstract}
Delsarte-Goethals-Seidel showed that if $X$ is a spherical $t$-design with degree $s$ satisfying $t\geq 2s-2$, 
$X$ carries the structure of an association scheme. 
Also Bannai-Bannai showed that the same conclusion holds if $X$ is an antipodal spherical $t$-design with degree $s$ satisfying $t=2s-3$.
As a generalization of these results, we prove that a union of spherical designs with a certain property carries 
the structure of a coherent configuration.
We derive triple regularity of tight spherical $4,5,7$-designs, mutually unbiased bases, 
linked symmetric designs with certain parameters. 
\end{abstract}

\section{Introduction}
Spherical codes and designs were studied by Delsarte-Goethals-Seidel \cite{DGS}.
There are two important parameters of finite set $X$ in the unit sphere $S^{d-1}$, that is, strength $t$ and degree $s$.  
In the paper \cite{DGS}, it is shown that $t\geq 2s-2$ implies $X$ carries an $s$-class association scheme.
Recently Bannai-Bannai \cite{BB} has shown that if $X$ is antipodal and $t=2s-3$, 
then $X$ carries an $s$-class association scheme.

Coherent configurations, that were introduced by D. G. Higman \cite{Hig}, are known as a generalization of association schemes.
In Section 2, as an analogue of these results, we give a certain sufficient condition for a union of spherical designs 
to carry the structure of a coherent configuration.
Our proof is based on the method of Delsarte-Goethals-Seidel~\cite[Theorem 7.4]{DGS}.

In Section 3, we consider triply regular association schemes which were introduced in connection with spin models by F. Jaeger \cite{J} and have higher regularity than ordinary association schemes.
Triple regularity is equivalent to the condition that
the partition consisting of subconstituents relative to any point of the association scheme carries a coherent configuration whose parameters are independent of the point.
In order to show that a symmetric association scheme is triply regular, we embed the scheme to the unit sphere $S^{d-1}$ by a primitive idempotent.
This embedding has a partition of derived designs in $S^{d-2}$ for arbitrary point in the association scheme.
Applying the main theorem of this paper to the union of derived designs, we obtain a sufficient condition for triple regularity of a symmetric association scheme.

In Sections \ref{sec:3}--\ref{sec:6}, we consider tight spherical $4,5,7$-designs, mutually unbiased bases (MUB), and linked symmetric designs with certain parameters. 
We note that tight spherical $t$-designs are classified except for $t=4,5,7$.
It is known that a tight spherical design, MUB,
and a linked system of symmetric designs carry a symmetric association scheme \cite[Theorem 7.4]{DGS}, \cite[Theorem 1.1]{BB}, \cite{M}.
We will show that these symmetric association schemes are triply regular using our main theorem.
\section{Coherent configurations obtained from spherical designs}\label{sec:2}
Let $X$ be a finite set, 
we define $\text{diag} (X\times X)=\{(x,x)\mid x\in X\}$.
Let $\{f_i\}_{i\in I}$ be a set of relations on $X$, 
we define $f_i^t=\{(y,x)\mid (x,y)\in f_i\}$.
$(X,\{f_i\}_{i\in I})$ is a coherent configuration if the following properties are satisfied:
\begin{enumerate}
\item $\{f_i\}_{i\in I}$ is a partition of $X\times X$,
\item $f_i^t=f_{i^*}$ for some $i^*\in I$,
\item $f_i \cap\mathrm{diag}(X\times X)\neq \emptyset$ implies $f_i \subset\mathrm{diag}(X\times X)$,
\item for $i,j,k\in I$, the number $|\{z\in X \mid (x,z)\in f_i,(z,y)\in f_j\}|$ is independent of the choice of $(x,y)\in f_k$.   
\end{enumerate}
If moreover $f_0=\mathrm{diag}(X\times X)$ and $i^*=i$ for all $i \in I$, then we call $(X,\{f_i\}_{i\in I})$ a symmetric association scheme.

Let $X_1,\dots,X_n$ be finite subsets of $S^{d-1}$.
We denote by $\coprod \nolimits_{i=1}^n X_i$ the disjoint union of $X_1,\dots,X_n$.
We denote by $\langle x,y\rangle$ the inner product of $x, y \in \mathbb{R}^d$.
We define the nontrivial angle set $A(X_i,X_j)$ between $X_i$ and $X_j$ by
$$A(X_i,X_j)=\{ \langle x, y\rangle \mid x \in X_i , y \in X_j , x \neq\pm y \},$$
and the angle set $A'(X_i,X_j)$ between $X_i$ and $X_j$ by
$$A'(X_i,X_j)=\{ \langle x, y\rangle \mid x \in X_i , y \in X_j, x\neq y\}.$$ 
If $i=j$, then $A(X_i,X_i)$ (resp. $A'(X_i,X_i)$) is abbreviated $A(X_i)$ (resp. $A'(X_i)$).

We define the intersection numbers on $X_j$ for $x,y\in S^{d-1}$ by
$$p_{\alpha,\beta}^j(x,y)=|\{z\in X_j \mid \langle x,z\rangle=\alpha,\langle y,z\rangle=\beta \}|.$$

For a positive integer $t$, a finite non-empty set $X$ in the unit sphere $S^{d-1}$ is called 
a spherical $t$-design in $S^{d-1}$ if the following condition is satisfied:
$$\frac{1}{|X|}\sum\limits_{x \in X}f(x)=\frac{1}{|S^{d-1}|}\int\nolimits_{S^{d-1}}f(x)d\sigma(x)$$
for all polynomials $f(x)=f(x_1,\dots,x_d)$ of degree not exceeding $t$.
Here $|S^{d-1}|$ denotes the volume of the sphere $S^{d-1}$.
When $X$ is a $t$-design and not a $(t+1)$-design, we call $t$ its strength. 

We define the Gegenbauer polynomials $\{Q_k(x)\}_{k=0}^\infty$ on $S^{d-1}$ by
\begin{align*}
& Q_0(x)=1,\quad Q_1(x)=dx,\\
& \frac{k+1}{d+2k}Q_{k+1}(x)=xQ_k(x)-\frac{d+k-3}{d+2k-4}Q_{k-1}(x).
\end{align*}

Let $\mbox{Harm}(\mathbb{R}^d)$ be the vector space of the harmonic polynomials over $\mathbb{R}$
and $\mbox{Harm}_l(\mathbb{R}^d)$ be the subspace of $\mbox{Harm}(\mathbb{R}^d)$ consisting of homogeneous polynomials of total degree $l$.
Let $\{\phi_{l,1},\dots,\phi_{l,h_l}\}$ be an orthonormal basis
of $\text{Harm}_l(\mathbb{R}^d)$ with respect to the inner product 
$$\langle\phi,\psi \rangle=\frac{1}{|S^{d-1}|}\int\nolimits_{S^{d-1}}\phi(x)\psi(x)d\sigma(x) .$$
Then the addition formula for the Gegenbauer polynomial holds \cite[Theorem 3.3]{DGS}:
\begin{lemma}\label{add}
$\sum\limits_{i=1}^{h_l}\phi_{l,i}(x)\phi_{l,i}(y)=Q_l(\langle x,y\rangle)$ for any $l\in \mathbb{N}$, $x,y\in S^{d-1}$.
\end{lemma}
We define the $l$-th characteristic matrix of a finite set $X\subset S^{d-1}$ as the $|X|\times h_l$ matrix
$$ H_l=(\phi_{l,i}(x))_{\substack{x\in X\\1\leq i\leq h_l}} .$$
A criterion for $t$-designs using Gegenbauer polynomials and the characteristic matrices is known \cite[Theorem 5.3, 5.5]{DGS}.
\begin{lemma}\label{cha}
Let $X$ be a finite set in $S^{d-1}$. The following conditions are equivalent:
\begin{enumerate}
\item $X$ is a $t$-design,
\item $\sum\limits_{x,y\in X}Q_k(\langle x,y\rangle)=0$ for any $k\in \{1,\ldots,t\}$,
\item $H_k^tH_l=\delta_{k,l}|X|I \quad \text{for} \quad 0\leq k+l\leq t $, 
\end{enumerate}
\end{lemma}

We define $\{f_{\lambda,l}\}_{l=0}^\lambda$ as the coefficients of Gegenbauer expansion of $x^\lambda$ for any nonnegative integers $\lambda$, i.e.,
$x^\lambda=\sum\nolimits_{l=0}^\lambda f_{\lambda,l}Q_l(x)$, 
and let
$F_{\lambda,\mu}(x)=\sum\nolimits_{l=0}^{\min\{\lambda,\mu\}}f_{\lambda,l}f_{\mu,l}Q_l(x)$, where $\lambda, \mu$ are nonnegative integers.

The following three lemmas are used to prove Theorem~\ref{coherent} by using uniqueness of the solution of linear equations.   
Let $A$ be a square matrix of size $n$. 
For index sets $I,J \subset\{1,\dots,n\}$, we denote the submatrix that lies in the rows of $A$ indexed by $I$ and the columns indexed by $J$ as $A(I,J)$ and the complement of $I$ as $I'$.
If $I=\{i\}$ and $J=\{j\}$, then $A(I,J)$ is abbreviated $A(i,j)$.
A lemma which relates a minor of $A^{-1}$ to that of $A$ is the following:   
\begin{lemma}\label{det}{\upshape \cite[p.21]{HJ}}
Let $A$ be a nonsingular matrix, 
and let $I,J$ be index sets of rows and columns of $A$ with $|I|=|J|$. Then
$$\det A^{-1}(I',J')=(-1)^{\sum\nolimits_{i\in I}i+\sum\nolimits_{j\in J}j}\frac{\det A(J,I)}{\det A}.$$
\end{lemma}
We define the $k$-th elementary symmetric polynomial $e_k(x_1,\ldots,x_n)$ in $n$ valuables $x_1,\ldots,x_n$ by  
$$e_k(x_1,\ldots,x_n)=
\begin{cases}1 & \text{if}\ k=0,\\
\sum\limits_{1\leq i_1<\cdots<i_k\leq n}
x_{i_1}x_{i_2}\cdots x_{i_k} & \text{if}\ k\geq1.
\end{cases}$$
We define the polynomial $a_{\lambda}(x_1,\ldots,x_n)$ for a partition $\lambda=(\lambda_1,\ldots,\lambda_n)$ by
$$a_{\lambda}(x_1,\ldots,x_n)=\sum\limits_{\sigma\in S_n}\epsilon(\sigma)x_{\sigma(1)}^{\lambda_1}\cdots x_{\sigma(n)}^{\lambda_n},$$
and the Schur function $S_{\lambda}(x_1,\ldots,x_n)$ by
$$S_{\lambda}(x_1,\ldots,x_n)=\frac{a_{\lambda+\delta}(x_1,\dots,x_n)}{a_{\lambda}(x_1,\dots,x_n)},$$
where $\delta=(n-1,n-2,\ldots,1,0)$.
\begin{lemma}\label{detv}
Let $A$ be a square matrix of order $n$ with $(i,j)$ entry
$\alpha_j^{i-1}$, 
where $\alpha_1,\cdots,\alpha_n$ are distinct. 
Then
$$A^{-1}(i,j)=(-1)^{i+j}\frac{e_{n-j}(\alpha_1,\ldots,\alpha_{i-1},\alpha_{i+1},\ldots,\alpha_n)}{\prod\limits_{1\leq k<i} (\alpha_i-\alpha_k)\prod\limits_{i<l\leq n}(\alpha_l-\alpha_i)}.$$
\end{lemma}
\begin{proof}
Putting $\lambda=(1^{n-j},0^{j-1})$, we have by \cite[p.42]{Mac}, 
\begin{align*}
A^{-1}(i,j)&=(-1)^{i+j}\frac{\det{A(\{j\}',\{i\}')}}{\det{A}} \\
&=(-1)^{i+j}\frac{a_{\lambda+\delta}(\alpha_1,\ldots,\alpha_{i-1},\alpha_{i+1},\ldots,\alpha_n)}{\det{A}}\\ \displaybreak[0]
&=\frac{(-1)^{i+j}}{\prod\limits_{1\leq k<i} (\alpha_i-\alpha_k)\prod\limits_{i<l\leq n}(\alpha_l-\alpha_i)}\frac{a_{\lambda+\delta}(\alpha_1,\ldots,\alpha_{i-1},\alpha_{i+1},\ldots,\alpha_n)}{a_{\delta}(\alpha_1,\ldots,\alpha_{i-1},\alpha_{i+1},\ldots,\alpha_n)}\\ \displaybreak[0]
&=\frac{(-1)^{i+j}}{\prod\limits_{1\leq k<i} (\alpha_i-\alpha_k)\prod\limits_{i<l\leq n}(\alpha_l-\alpha_i)}S_{\lambda}(\alpha_1,\ldots,\alpha_{i-1},\alpha_{i+1},\ldots,\alpha_n)\\ \displaybreak[0]
&=\frac{(-1)^{i+j}}{\prod\limits_{1\leq k<i} (\alpha_i-\alpha_k)\prod\limits_{i<l\leq n}(\alpha_l-\alpha_i)}e_{n-j}(\alpha_1,\ldots,\alpha_{i-1},\alpha_{i+1},\ldots,\alpha_n) \displaybreak[0]
\end{align*}
\end{proof}
\begin{lemma}\label{detv3}
Let $A$ be a square matrix of order $n$ with $(i,j)$ entry
$\alpha_j^{i-1}$ 
and
Let $B$ be a square matrix of order $m$ with $(i,j)$ entry
$\beta_j^{i-1}$, where
$\alpha_1,\cdots,\alpha_n$ and
$\beta_1,\cdots,\beta_m$ are distinct.
Let $J,I$ be index sets of rows and columns, respectively, of
$A\otimes B$ such that
$J'=\{(n-1,m),(n,m-1),(n,m)\}$,
$I'=\{(i_1,j_1),(i_2,j_2),(i_3,j_3)\}$.
Then
$$\frac{\det{(A\otimes B)(J,I)}}{\det{A\otimes B}}=\pm\frac{\alpha_{i_1}\beta_{j_2}+\alpha_{i_2}\beta_{j_3}+\alpha_{i_3}\beta_{j_1}-\alpha_{i_1}\beta_{j_3}-\alpha_{i_2}\beta_{j_1}-\alpha_{i_3}\beta_{j_2}}{\prod\limits_{1\leq r\leq 3}\biggl(\prod\limits_{1\leq k<i_r}(\alpha_{i_r}-\alpha_k)\prod\limits_{i_r<l\leq n}(\alpha_l-\alpha_{i_r})\prod\limits_{1\leq k<j_r} (\beta_{j_r}-\beta_k)\prod\limits_{j_r<l\leq m}(\beta_l-\beta_{j_r})\biggr)}.$$
\end{lemma}
\begin{proof}
We define $f(i,j)=\prod\limits_{1\leq k<i}(\alpha_i-\alpha_k)\prod\limits_{i<l\leq n}(\alpha_l-\alpha_i)\prod\limits_{1\leq k<j} (\beta_j-\beta_k)\prod\limits_{j<l\leq m}(\beta_l-\beta_j)$.
Using Lemmas~\ref{det} and \ref{detv},
\begin{align*}
\frac{\det{(A\otimes B)(J,I)}}{\det{A\otimes B}}&=\pm\det{(A\otimes B)^{-1}(I',J')}\\
&=\pm\det{(A^{-1}\otimes B^{-1})(I',J')}\displaybreak[0]\\
&=\pm\det{\left( \begin{array}{ccc}
        \frac{(-1)^{i_1+n-1+j_1+m}\sum\nolimits_{i\neq i_1}\alpha_i}{f(i_1,j_1)}  & \frac{(-1)^{i_1+n+j_1+m-1}\sum\nolimits_{j\neq j_1}\beta_j}{f(i_1,j_1)} & \frac{(-1)^{i_1+n+j_1+m}}{f(i_1,j_1)} \\
        \frac{(-1)^{i_2+n-1+j_2+m}\sum\nolimits_{i\neq i_2}\alpha_i}{f(i_2,j_2)}  & \frac{(-1)^{i_2+n+j_2+m-1}\sum\nolimits_{j\neq j_2}\beta_j}{f(i_2,j_2)} & \frac{(-1)^{i_2+n+j_2+m}}{f(i_2,j_2)} \\
        \frac{(-1)^{i_3+n-1+j_3+m}\sum\nolimits_{i\neq i_3}\alpha_i}{f(i_3,j_3)}  & \frac{(-1)^{i_3+n+j_3+m-1}\sum\nolimits_{j\neq j_3}\beta_j}{f(i_3,j_3)} & \frac{(-1)^{i_3+n+j_3+m}}{f(i_3,j_3)}   \end{array}\right)} \displaybreak[0]\\
&=\pm\frac{1}{\prod\limits_{1\leq r\leq 3} f(i_r,j_r)}\det{\left(\begin{array}{ccc} \sum\nolimits_{i\neq i_1}\alpha_i&\sum\nolimits_{j\neq j_1}\beta_j&1\\ \sum\nolimits_{i\neq i_2}\alpha_i&\sum\nolimits_{j\neq j_2}\beta_j&1\\ \sum\nolimits_{i\neq i_3}\alpha_i&\sum\nolimits_{j\neq j_3}\beta_j&1 \end{array}\right)}     \displaybreak[0] \\
&=\pm\frac{1}{\prod\limits_{1\leq r\leq 3} f(i_r,j_r)}\det{\left(\begin{array}{ccc} \alpha_{i_1}&\beta_{j_1}&1\\ \alpha_{i_2}&\beta_{j_2}&1\\ \alpha_{i_3}&\beta_{j_3}&1 \end{array}\right)}  \displaybreak[0]\\
&=\pm\frac{\alpha_{i_1}\beta_{j_2}+\alpha_{i_2}\beta_{j_3}+\alpha_{i_3}\beta_{j_1}-\alpha_{i_1}\beta_{j_3}-\alpha_{i_2}\beta_{j_1}-\alpha_{i_3}\beta_{j_2}}{\prod\limits_{1\leq r\leq 3}\biggl(\prod\limits_{1\leq k<i_r}(\alpha_{i_r}-\alpha_k)\prod\limits_{i_r<l\leq n}(\alpha_l-\alpha_{i_r})\prod\limits_{1\leq k<j_r} (\beta_{j_r}-\beta_k)\prod\limits_{j_r<l\leq m}(\beta_l-\beta_{j_r})\biggr)}.
\end{align*}
\end{proof}
The following is the main theorem of this paper.
\begin{theorem}\label{coherent}
Let $X_i \subset S^{d-1}$ be a spherical $t_i$-design for $i \in \{1,\dots,n\}$.
Assume that $X_i \cap X_j=\emptyset$ or $X_i=X_j$, and $X_i \cap (-X_j)=\emptyset$ or $X_i=-X_j$ for $i ,j \in \{1,\dots,n\}$.
Let $s_{i,j}=|A(X_i,X_j)|$, $s_{i,j}^*=|A'(X_i,X_j)|$ and $A(X_i,X_j)=\{\alpha_{i,j}^1,\ldots,\alpha_{i,j}^{s_{i,j}}\}$, $\alpha_{i,j}^0=1$, when $-1\in A'(X_i,X_j)$, we define  $\alpha_{i,j}^{s_{i,j}^*}=-1$. 
We define $R_{i,j}^k=\{(x,y)\in X_i \times X_j\mid \langle x,y \rangle =\alpha_{i,j}^k \} $.
If one of the following holds depending on the choice of $i,j,k \in \{1, \dots ,n \}$:
\begin{enumerate}
\item $s_{i,j}+s_{j,k}-2 \leq t_j$,
\item $s_{i,j}+s_{j,k}-3 = t_j$ and for any $\gamma \in A(X_i,X_k)$ there exist $\alpha\in A(X_i,X_j),\beta\in A(X_j,X_k)$ such that the number $p_{\alpha,\beta}^j(x,y)$ is independent of the choice of $x\in X_i,y\in X_k$ with $\gamma=\langle x,y\rangle$,
\item $s_{i,j}+s_{j,k}-4 = t_j$ and for any $\gamma \in A(X_i,X_k)$  there exist $\alpha,\alpha'\in A(X_i,X_j),\beta,\beta'\in A(X_j,X_k)$  such that $\alpha\neq \alpha'$, $\beta\neq \beta'$ and 
the numbers $p_{\alpha,\beta}^j(x,y)$, $p_{\alpha,\beta'}^j(x,y)$ and $p_{\alpha',\beta}^j(x,y)$ are independent of the choice of $x\in X_i,y\in X_k$ with $\gamma=\langle x,y\rangle$,
\end{enumerate}   
then $(\coprod \nolimits_{i=1}^n X_i, \{R_{i,j}^k\mid 1\leq i,j\leq n,1-\delta_{X_i,X_j}\leq k\leq s_{i,j}^*\}) $ is a coherent configuration.
The parameters of this coherent configuration are determined by $A(X_i,X_j)$, $|X_i|$, $t_i$, $\delta_{X_i,X_j}$, $\delta_{X_i,-X_j}$, and when $s_{i,j}+s_{j,k}-3=t_j$ (resp. $s_{i,j}+s_{j,k}-4=t_j$), the numbers $p_{\alpha,\beta}^j(x,y)$ (resp. $p_{\alpha,\beta}^j(x,y)$, $p_{\alpha',\beta}^j(x,y)$, $p_{\alpha,\beta'}^j(x,y)$) which are assumed be independent of $(x,y)$ with $\langle x,y\rangle=\gamma$. 
\end{theorem}
\begin{proof}
Let $x\in X_i$, $y\in X_k$ be such that $\gamma=\langle x,y\rangle$. 
It is sufficient to show that the number $p_{\alpha,\beta}^j(x,y)$ depends only on $\gamma$ and does not depend on the choice of  $x\in X_i,y\in X_k$ satisfying $\gamma=\langle x,y\rangle$.

For the ease of notation, let $\alpha_l=\alpha_{i,j}^l$ and $\beta_m=\alpha_{j,k}^m$.

We define a mapping $\phi_l:S^{d-1}\rightarrow\mathbb{R}^{h_l}$ by $\phi_l(x)=(\varphi_{l,1}(x),\dots,\varphi_{l,h_l}(x))$.
Let $H_l$ be the $l$-th characteristic matrix of $X_j$. 
For any non-negative integers $\lambda$ and $\mu$ satisfying $\lambda +\mu \leq t_j$, we calculate 
$$(\sum \limits_{l=1}^\lambda f_{\lambda,l}\phi_l(x)H_l^t)(\sum \limits_{m=1}^\mu f_{\mu,m}H_m\phi_m(y)^t)$$
in two different ways.

First we use Lemma~\ref{cha} and Lemma~\ref{add} in turn, to obtain the following equality:
\begin{align}
(\sum \limits_{l=1}^\lambda f_{\lambda,l}\phi_l(x)H_l^t)(\sum \limits_{m=1}^\mu f_{\mu,m}H_m\phi_m(y)^t) &=|X_j|\sum \limits_{l=1}^{\min\{\lambda,\mu\}}  f_{\lambda,l} f_{\mu,l}\phi_l(x)\phi_l(y)^t\notag \\
&=|X_j|\sum \limits_{l=1}^{\min\{\lambda,\mu\}}  f_{\lambda,l} f_{\mu,l}Q_l(\langle x,y\rangle)\notag \\
&=|X_j|F_{\lambda,\mu}(\langle x,y\rangle) \label{eq1}. 
\end{align}
Next using Lemma~\ref{add}, we obtain the following equality:
\begin{align}
\lefteqn{(\sum \limits_{l=1}^\lambda f_{\lambda,l}\phi_l(x)H_l^t)(\sum \limits_{m=1}^\mu f_{\mu,m}H_m\phi_m(y)^t)} \notag \\
&= \sum \limits_{z \in X_j}(\sum \limits_{l=1}^\lambda f_{\lambda,l}(\phi_l(x)\phi_l(z)^t)(\sum \limits_{m=1}^\mu f_{\mu,m}(\phi_m(z) \phi_m(y)^t)\notag \displaybreak[0]\\
&=\sum \limits_{z \in X_j}(\sum \limits_{l=1}^\lambda f_{\lambda,l}Q_l(\langle x,z \rangle))(\sum \limits_{m=1}^\mu f_{\mu,m}Q_m(\langle z,y \rangle)) \notag \displaybreak[0]\\
&=\sum \limits_{z \in X_j}\langle x,z \rangle^\lambda \langle z,y \rangle^\mu \notag \displaybreak[0]\\
&=\sum\limits_{\substack{\alpha \in A'(X_i,X_j)\\ \beta \in A'(X_j,X_k)}}\alpha^\lambda \beta^\mu p_{\alpha,\beta}^j(x,y)+p_{1,1}^j(x,y)+ \sum\limits_{m=1}^{s_{j,k}^*} \beta_m^\mu p_{1,\beta_m}^j(x,y)+\sum\limits_{l=1}^{s_{i,j}^*}\alpha_l^\lambda  p_{\alpha_l,1}^j(x,y) \notag \displaybreak[0]\\ 
&=\sum\limits_{l=1}^{s_{i,j}}\sum\limits_{m=1}^{s_{j,k}}\alpha_l^\lambda \beta_m^\mu p_{\alpha_l,\beta_m}^j(x,y) \notag \\
&\quad+ p_{1,1}^j(x,y)+ (-1)^\mu p_{1,-1}^j(x,y)+(-1)^\lambda p_{-1,1}^j(x,y)+(-1)^\lambda (-1)^\mu p_{-1,-1}^j(x,y) \notag \\
&\quad+\sum\limits_{m=1}^{s_{j,k}} \beta_m^\mu p_{1,\beta_m}^j(x,y)+\sum\limits_{l=1}^{s_{i,j}}\alpha_l^\lambda  p_{\alpha_l,1}^j(x,y)+ \sum\limits_{m=1}^{s_{j,k}}(-1)^\lambda \beta_m^\mu p_{-1,\beta_m}^j(x,y)+\sum\limits_{l=1}^{s_{i,j}}\alpha_l^\lambda (-1)^\mu p_{\alpha_l,-1}^j(x,y) \notag \displaybreak[0]\\
&=\sum\limits_{l=1}^{s_{i,j}}\sum\limits_{m=1}^{s_{j,k}}\alpha_l^\lambda \beta_m^\mu p_{\alpha_l,\beta_m}^j(x,y)+G_{\lambda,\mu}^{i,j,k}(\gamma),  \label{eq2}
\end{align}
where 
\begin{align*}
G_{\lambda,\mu}^{i,j,k}(t)&= \delta_{1,t}\delta_{X_i,X_j}\delta_{X_j,X_k}+(-1)^\mu\delta_{-1,t}\delta_{X_i,X_j}\delta_{X_j,-X_k}\\
&+(-1)^\lambda\delta_{-1,t}\delta_{X_i,-X_j}\delta_{X_j,X_k}+(-1)^{\lambda+\mu} \delta_{1,t}\delta_{X_i,-X_j}\delta_{X_j,-X_k}\\
&+(1-\delta_{1,t})(1-\delta_{-1,t})(\delta_{X_i,X_j}t^\mu +\delta_{X_j,X_k}t^\lambda+\delta_{X_i,-X_j}(-1)^\lambda (-t)^\mu+\delta_{X_j,-X_k}(-t)^\lambda (-1)^\mu).
\end{align*}
We obtain from \eqref{eq1} and \eqref{eq2}:
\begin{align} \sum\limits_{l=1}^{s_{i,j}}\sum\limits_{m=1}^{s_{j,k}}\alpha_l^\lambda \beta_m^\mu p_{\alpha_l,\beta_m}^j(x,y)  &=
   |X_j|F_{\lambda,\mu}(\langle x,y\rangle)-G_{\lambda,\mu}^{i,j,k}(\langle x,y\rangle).\label{eq3}
\end{align}

In the case where $i,j,k$ satisfy the assumption $(1)$, for $0 \leq \lambda \leq s_{i,j}-1$ and $0 \leq \mu \leq s_{j,k}-1$, (\ref{eq3}) yields a system of $s_{i,j}s_{j,k}$ linear equations whose unknowns are $$\{p_{\alpha_l,\beta_m}^j(x,y)\mid 1\leq l\leq s_{i,j},\ 1\leq m\leq s_{j,k}\}.$$
Its coefficient matrix  $A\otimes B$ is nonsingular, where 
$$A=
 \left( \begin{array}{ccc} 
 1     &  \cdots & 1 \\
 \alpha_1  & \cdots & \alpha_{s_{i,j}} \\ 
 \vdots & \ddots & \vdots \\
 \alpha_1^{s_{i,j}-1} & \cdots & \alpha_{s_{i,j}}^{s_{i,j}-1}
\end{array} \right),\quad B=\left( \begin{array}{ccc} 
 1     &  \cdots & 1 \\
 \beta_1  & \cdots & \beta_{s_{j,k}} \\ 
 \vdots & \ddots & \vdots \\
 \beta_1^{s_{j,k}-1} & \cdots & \beta_{s_{j,k}}^{s_{j,k}-1}
\end{array} \right).$$  
Therefore $p_{\alpha_l,\beta_m}^j(x,y)$ for $ 1\leq l\leq s_{i,j}$, $1\leq m\leq s_{j,k}$ depends only on $\gamma$ and does not depend on the choice of $x$, $y$ satisfying $\gamma=\langle x,y\rangle$, 
and is determined by $A(X_i,X_j)$, $A(X_j,X_k)$, $\gamma$, $|X_j|$, $t_j$, $\delta_{X_i,X_j}$, $\delta_{X_j,X_k}$, $\delta_{X_i,-X_j}$, $\delta_{X_j,-X_k}$.

In the case where $i,j,k$ satisfy $(2)$ i.e., for  $\langle x,y \rangle =\gamma \in A(X_i,X_k)$, there exist $\alpha_{l^*}\in A(X_i,X_j)$, $\beta_{m^*}\in A(X_j,X_k)$ such that
the number $p_{\alpha_{l^*},\beta_{m^*}}^j(x,y)$ is uniquely determined. 
The linear equation~\eqref{eq3} is the following:
\begin{align} \sum\limits_{\substack{1 \leq l \leq s_{i,j}\\ 1 \leq m \leq s_{j,k}\\ (l,m)\neq (l^*,m^*)}}\alpha_l^\lambda \beta_m^\mu p_{\alpha_l,\beta_m}^j(x,y)  &=
   |X_j|F_{\lambda,\mu}(\langle x,y\rangle)-G_{\lambda,\mu}^{i,j,k}(\langle x,y\rangle)-\alpha_{l^*}^\lambda \beta_{m^*}^\mu p_{\alpha_{l^*},\beta_{m^*}}^j(x,y). \label{eq4}
\end{align}
For $0 \leq \lambda \leq s_{i,j}-1, 0 \leq \mu \leq s_{j,k}-1$ and $(\lambda,\mu) \neq (s_{i,j}-1,s_{j,k}-1) $, \eqref{eq4} yields a system of $s_{i,j}s_{j,k}-1$ linear equations whose unknowns are $$\{p_{\alpha_l,\beta_m}^j(x,y)\mid
1\leq l\leq s_{i,j},\; 1\leq m\leq s_{j,k},\; (l,m)\neq(l^*,m^*)\}.$$
The coefficient matrix $C_1$ of these linear equations is the submatrix obtained by deleting the $(s_{i,j},s_{j,k})$-row and $(l^*,m^*)$-column of $A\otimes B$. 
Using Lemma~\ref{detv} the determinant of $C_1$ is, up to sign, 
\begin{align*}
\det{C_1}&= \pm ((s_{i,j},s_{j,k}),(l^*,m^*))\mbox{-cofactor of}\ A\otimes B   \\
          &=\pm ((l^*,m^*),((s_{i,j},s_{j,k}))\mbox{-entry of}\ (A\otimes B)^{-1})\det{A\otimes B} \displaybreak[0]\\ 
          &=\pm ((l^*,s_{i,j})\mbox{-entry of}\ A^{-1})\times ((m^*,s_{j,k})\mbox{-entry of}\ B^{-1})\det{A\otimes B} \displaybreak[0]\\ 
          &=\pm\frac{\det{A\otimes B}}{\prod\limits_{1\leq k<l^*} (\alpha_{l^*}-\alpha_k) \prod\limits_{l^*<l\leq s_{i,j}} (\alpha_l-\alpha_{l^*})\prod\limits_{1\leq k<m^*} (\beta_{m^*}-\beta_k) \prod\limits_{m^*<l\leq s_{j,k}} (\beta_l-\beta_{m^*})}.
\end{align*}          
Hence $C_1$ is nonsingular.

Therefore $p_{\alpha_l,\beta_m}^j(x,y)$ for $ 1\leq l\leq s_{i,j}$, $1\leq m\leq s_{j,k}$, $(l,m)\neq(l^*,m^*)$ depends only on $\gamma$ and does not depend on the choice of $x$, $y$ satisfying $\gamma=\langle x,y\rangle$, 
and is determined by $A(X_i,X_j)$, $A(X_j,X_k)$, $\gamma$, $|X_j|$, $t_j$, $\delta_{X_i,X_j}$, $\delta_{X_j,X_k}$, $\delta_{X_i,-X_j}$, $\delta_{X_j,-X_k}$, the number $p_{\alpha_{l^*},\beta_{m^*}}^j(x,y)$
which is assumed be independent of $(x,y)$ with $\langle x,y\rangle=\gamma$.\\

In the case where $i,j,k$ satisfy $(3)$ i.e., for  $\langle x,y \rangle =\gamma \in A(X_i,X_k)$ there exist $\alpha_{l_1},\alpha_{l_2}\in A(X_i,X_j)$, $\beta_{m_1},\beta_{m_2}\in A(X_j,X_k)$ such that
the numbers $p_{\alpha_{l_1},\beta_{m_1}}^j(x,y),p_{\alpha_{l_1},\beta_{m_2}}^j(x,y),p_{\alpha_{l_2},\beta_{m_1}}^j(x,y)$ are uniquely determined. 
The linear equation~\eqref{eq3} is the following:
\begin{align} \sum\limits_{\substack{1 \leq l \leq s_{i,j}\\ 1 \leq m \leq s_{j,k}\\ (l,m)\neq (l_1,m_1),(l_1,m_2),(l_2,m_1)}}\!\!\!\!\!\!\!\!\!\!\!\!\!\!\!\!\!\!\!\!\!\!\!\!\alpha_l^\lambda \beta_m^\mu p_{\alpha_l,\beta_m}^j(x,y)  &=
   |X_j|F_{\lambda,\mu}(\langle x,y\rangle)-G_{\lambda,\mu}^{i,j,k}(\langle x,y\rangle)-\alpha_{l_1}^\lambda \beta_{m_1}^\mu p_{\alpha_{l_1},\beta_{m_1}}^j(x,y)\notag \\[-10mm] 
&\quad-\alpha_{l_1}^\lambda \beta_{m_2}^\mu p_{\alpha_{l_1},\beta_{m_2}}^j(x,y)-\alpha_{l_2}^\lambda \beta_{m_1}^\mu p_{\alpha_{l_2},\beta_{m_1}}^j(x,y). \label{eq5}
\end{align}
For $0 \leq \lambda \leq s_{i,j}-1, 0 \leq \mu \leq s_{j,k}-1$ and $(\lambda,\mu)\neq (s_{i,j}-2,s_{j,k}-1),(s_{i,j}-1,s_{j,k}-2),(s_{i,j}-1,s_{j,k}-1)$, 
\eqref{eq5} yields a system of $s_{i,j}s_{j,k}-3$ linear equations whose unknowns are 
$$\{p_{\alpha_l,\beta_m}^j(x,y)\mid
1\leq l\leq s_{i,j},\; 1\leq m\leq s_{j,k},\; (l,m)\neq(l_1,m_1),(l_1,m_2),(l_2,m_1)\}.$$
The coefficient matrix $C_2$ of these linear equations is the submatrix obtained by deleting the $(s_{i,j}-1,s_{j,k}),(s_{i,j},s_{j,k}-1 ),(s_{i,j},s_{j,k})$-rows and $(l_1,m_1),(l_1,m_2),(l_2,m_1)$-columns of $A\otimes B$.
Let $J,I$ be index sets of rows and columns, respectively, of
$A\otimes B$ such that
$$J'=\{(s_{i,j}-1,s_{j,k}),(s_{i,j},s_{j,k}-1 ),(s_{i,j},s_{j,k})\}$$ and $$I'=\{(l_1,m_1),(l_1,m_2),(l_2,m_1)\}.$$
Setting $(i_1,j_1),(i_2,j_2),(i_3,j_3)$ to be
$(l_1,m_1),(l_1,m_2),(l_2,m_1)$ respectively,
we have
$$
\alpha_{i_1}\beta_{j_2}+\alpha_{i_2}\beta_{j_3}+\alpha_{i_3}\beta_{j_1}
-\alpha_{i_1}\beta_{j_3}-\alpha_{i_2}\beta_{j_1}-\alpha_{i_3}\beta_{j_2}
=(\alpha_{l_1}-\alpha_{l_2})(\beta_{m_1}-\beta_{m_2}).
$$
Hence $C_2$ is nonsingular by Lemma~\ref{detv3}.
Therefore $p_{\alpha_l,\beta_m}^j(x,y)$ for $ 1\leq l\leq s_{i,j}$, $1\leq m\leq s_{j,k}$, $(l,m)\neq(l_1,m_1),(l_1,m_2),(l_2,m_1)$ depends only on $\gamma$ and does not depend on the choice of $x$, $y$ satisfying $\gamma=\langle x,y\rangle$, 
and is determined by $A(X_i,X_j)$, $A(X_j,X_k)$, $\gamma$, $|X_j|$, $t_j$, $\delta_{X_i,X_j}$, $\delta_{X_j,X_k}$, $\delta_{X_i,-X_j}$, $\delta_{X_j,-X_k}$,
the numbers $p_{\alpha,\beta}^j(x,y)$, $p_{\alpha',\beta}^j(x,y)$, $p_{\alpha,\beta'}^j(x,y)$ which are assumed be independent of $(x,y)$ with $\langle x,y\rangle=\gamma$.
\end{proof}

Several results known for the case $n=1$ are derived from  Theorem~\ref{coherent}.
We consider the case where $n=1$ and $X=X_1$ is a $t$-design of degree $s$. 
Then $t_1=t$ and $$s_{1,1}=
\begin{cases}s-1 & \text{if}\ X \text{ is antipodal},\\
s & \text{if}\ X \text{ otherwise}.
\end{cases}$$
Suppose $t\geq 2s-2$.
If $X$ is antipodal, then $t_1\geq 2s_{1,1}$, 
and if $X$ is not antipodal, then $t_1\geq 2s_{1,1}-2$. 
Thus $X$ satisfies the assumption (1) of Theorem~\ref{coherent},
and hence $X$ carries a symmetric association scheme.
So Theorem~\ref{coherent} contains the first half of \cite[Theorem 7.4]{DGS} as a special case.

Suppose $t=2s-3$ and $p_{\gamma,\gamma}(x,y)$ is uniquely determined for any fixed $\gamma=\langle x,y\rangle\in A'(X)$.
If $X$ is antipodal, then $t_1=2s_{1,1}-1$, 
and if $X$ is not antipodal, then $t_1=2s_{1,1}-3$. 
Thus $X$ also satisfies the assumption (1) or (2) of
Theorem~\ref{coherent}, 
and hence $X$ carries a symmetric association scheme.
So Theorem~\ref{coherent} contains the second half of \cite[Theorem 7.4]{DGS} as a special case.

Suppose that $t=2s-3$.
If $X$ is antipodal, then $t_1=2s_{1,1}-1$.
Thus $X$ satisfies the assumption (1) of Theorem~\ref{coherent},
and hence $X$ carries a symmetric association scheme.
So Theorem~\ref{coherent} contains \cite[Theorem 1.1]{BB} as a special case.

Next, we consider triple regularity of a symmetric association scheme.
This concept was introduced in connection with spin models \cite{J}.

\begin{definition}
Let $(X,\{R_i\}_{i=0}^d)$ be a symmetric association scheme. Then the association scheme $X$ is said to be triply regular if,  
for all $i,j,k,l,m,n \in \{0,1,\dots,d \}$, and for all $x,y,z \in X$ such that $(x,y)\in R_i, (y,z) \in R_j, (z,x)\in R_k$, the number $p_{l,m,n}^{i,j,k}:=\lvert\{w \in X \mid (w,x) \in R_m, (w,y) \in R_n, (w,z) \in R_l \}\rvert$ depends only on $i,j,k,l,m,n$ and not on $x,y,z$.
\end{definition}
Let $(X,\{R_i\}_{i=0}^d)$ be an  association scheme. 
We define the $i$-th subconstituent with respect to $z\in X$ by $R_i(z):=\{y\in X \mid (z,y)\in R_i\}$. 
We denote by $R_{i,j}^k(z)$ the restriction of $R_k$ to $R_i(z)\times R_j(z)$.
The following lemma gives an equivalent definition of a triply regular association scheme.
We omit its easy proof.
\begin{lemma}\label{symtri}
A symmetric association scheme $(X,\{R_i\}_{i=0}^d)$
is triply regular if and only if
 for all $z\in X$, 
$(\bigcup\nolimits_{i=1}^dR_i(z), \{R_{i,j}^k(z)\mid 1\leq i,j\leq d,0\leq k \leq d,p_{i,j}^k\neq0\})$ 
is a coherent configuration whose parameters are independent of $z$. 
\end{lemma}
Let $X$ be a spherical $t$-design in $S^{d-1}$ with degree $s$, 
and $A'(X)=\{\alpha_1,\ldots,\alpha_s\}$. 
For $z\in X$ and $i\in\{1,\dots,s\}$,
$X_i(z)$ will denote the orthogonal projection of $\{y\in X \mid \langle y,z\rangle =\alpha_i\}$ to 
$z^\perp=\{y\in \mathbb{R}^d \mid \langle y,z\rangle =0 \}$,
rescaled to lie in $S^{d-2}$ in $z^\perp$.
$X_i(z)$ is called the derived design. 
In fact $X_i(z)$ is a $(t+1-s^*)$-design by \cite[Theorem 8.2]{DGS}, where $s^*=\lvert A'(X)\setminus \{-1\} \rvert$.
We define $\alpha_{i,j}^k=\frac{\alpha_k-\alpha_i\alpha_j}{\sqrt{(1-\alpha_i^2)(1-\alpha_j^2)}}$.
If $\langle x,z\rangle=\alpha_i$, $\langle y,z\rangle=\alpha_j$
and $\langle x,y\rangle=\alpha_k$,
then the inner product of the orthogonal projection of $x,y$ to
$z^\perp$ rescaled to lie in $S^{d-2}$, is $\alpha_{i,j}^k$. 
\begin{corollary}\label{triple}
Let $X\subset S^{d-1}$ be a finite set and $A'(X)=\{ \alpha_1,\dots,\alpha_s\}$.
Assume that $(X,\{R_k\}_{k=0}^s)$ is a symmetric association scheme,
where $R_k=\{(x,y)\in X\times X \mid \langle x,y\rangle=\alpha_k\}$
$(0\leq k\leq s)$ and $\alpha_0=1$.
Then 
\begin{enumerate}
\item $A(X_i(z),X_j(z))=\{\alpha_{i,j}^k \mid 0\leq k\leq s, p_{i,j}^k\neq0, \alpha_{i,j}^k\neq\pm1\}$.
\item $X_i(z)=X_j(z)$ or $X_i(z)\cap X_j(z)=\emptyset$,
and $X_i(z)=-X_j(z)$ or $X_i(z)\cap -X_j(z)=\emptyset$
for any $z\in X$ and any $i,j\in \{1,\ldots,s\}$. And $\delta_{X_i(z),X_j(z)}$, $\delta_{X_i(z),-X_j(z)}$ are independent of $z\in X$.
\item $X_i(z)$ has the same strength for all $z\in X$.
\end{enumerate}
Moreover if the assumption $(1),(2)$ or $(3)$ of Theorem~\ref{coherent}
is satisfied for $\{X_i(z)\}_{i=1}^s$,
and when $(i,j,k)$ satisfies $(2)$ (resp. $(3)$) the numbers $p_{\alpha,\beta}^j(x,y)$ (resp. $p_{\alpha,\beta}^j(x,y),p_{\alpha,\beta'}^j(x,y),p_{\alpha',\beta}^j(x,y)$) which are assumed to be independent of $(x,y)$ with $\gamma=\langle x,y\rangle$ are independent of the choice of $z$, 
then $(X,\{R_k\}_{k=0}^s)$ is a triply regular association scheme.
\end{corollary}
\begin{proof}
Let $z\in X$. (1) is immediate from the definition of $\alpha_{i,j}^k$. 

We define $R_{i,j}^k(z)=\{(x,y)\in X_i(z)\times X_j(z)\mid \langle x,y\rangle=\alpha_{i,j}^k\}$.
Then
\begin{align*}
\{\langle x,y\rangle\mid x\in &X_i(z),y\in X_j(z)\}\ni \pm 1 \\
\Leftrightarrow &\exists k\ \alpha_{i,j}^k=\pm1 \text{ and } p_{i,j}^k\neq 0   \\
\Leftrightarrow &\exists k\ \alpha_{i,j}^k=\pm 1 \text{, and}\\
&\qquad
\forall x\in X_i(z)\ \exists y\in X_j(z) \text{ s.t. } (x,y)\in R_{i,j}^k(z) \text{ and }
\\ &\qquad
\forall y\in X_j(z)\ \exists x\in X_i(z) \text{ s.t. } (x,y)\in R_{i,j}^k(z) \\
\Leftrightarrow &X_i(z)=\pm X_j(z).
\end{align*}
Since
$$\{\langle x,y\rangle\mid x\in X_i(z),y\in X_j(z)\}
=\{\alpha_{i,j}^k\mid 0\leq k\leq s, p_{i,j}^k\neq0\}$$ 
is independent of $z\in X$, $(2)$ holds.

By Lemma~\ref{cha}, $X_i(z)$ is a spherical $t$-design if and only if $\sum\nolimits_{x,y\in X_i(z)}Q_k(\langle x,y\rangle)=0$ for $k=1,\ldots,t$.
Since the number of $y\in X_i(z)$ satisfying $\langle x,y\rangle =\frac{\alpha_j-\alpha_i^2}{1-\alpha_i^2}$  is $p_{i,j}^i$
for any $x\in X_i(z)$,
the latter condition is equivalent to $\sum\nolimits_{0\leq j\leq s}Q_k(\frac{\alpha_j-\alpha_i^2}{1-\alpha_i^2})p_{i,j}^i=0$ for $k=1,\ldots,t$, which is independent of $z$.
Hence $X_i(z)$ has the same strength for all $z\in X$. 
Therefore $(3)$ holds.

Moreover if the assumption $(1),(2)$ or $(3)$ of Theorem~\ref{coherent}
is satisfied for $\{X_i(z)\}_{i=1}^s$, 
then $(\coprod\nolimits_{i=1}^sX_i(z), \{R_{i,j}^k(z)\mid 0\leq i,j,k\leq s,p_{i,j}^k\neq0\})$ is a coherent configuration. 
Clearly, $|X_i(z)|$ is  independent of $z\in X$.
Also, $A(X_i(z),X_j(z))$ is independent of $z\in X$ by (1),
$t_i$ is independent of $z\in X$ by (3), and
$\delta_{X_i(z),X_j(z)}$, $\delta_{X_i(z),-X_j(z)}$
are independent of $z\in X$ by (2).
It follows from Theorem~\ref{coherent} that the parameters of the coherent configuration are independent of $z\in X$.
Therefore, $(X,\{R_k\}_{k=0}^s)$ is a triply regular association scheme
by Lemma~\ref{symtri}.
\end{proof}

\section{Tight designs}\label{sec:3}
Let $X$ be a $t$-design in $S^{d-1}$.
It is known \cite[Theorems 5.11, 5.12]{DGS} that there is a lower bound  for the size of a spherical $t$-design in $S^{d-1}$.
Namely, if $X$ is a spherical $t$-design, then
$$|X|\geq \binom{d+t/2-1}{t/2} +\binom{n+t/2-2}{t/2-1}$$
if $t$ is even, and
$$|X|\geq 2\binom{d+(t-3)/2}{(t-1)/2} $$
if $t$ is odd.
If $X$ is a $t$-design for which one of the lower bounds is attained, then $X$ is called a tight $t$-design.
It was proved in \cite{BD1,BD2,DGS} that if $X$ is a tight $t$-design with degree $s$ in $S^{d-1}$,
then the following statements hold.
\begin{enumerate}
\item if $t$ is even, then $t=2s$,
\item if $t$ is odd, then $t=2s-1$ and $X$ is antipodal,
\item if $d=2$, then $X$ is the regular $(t+1)$-gon,
\item if $d \geq 3 $, then $t \leq 5$ or $t=7$, $11$.
\end{enumerate}
If $X$ is a tight $11$-design in $S^{d-1}$ where $d\geq3$, then $d=24$ and $X$ is the set of minimum vectors of
the Leech lattice \cite{BS}.  
We consider  tight $4$-, $5$-, $7$-designs in $S^{d-1}$ where $d\geq3$.

Let $X \subset S^{d-1}$ be a tight $2s$-design, and let 
$A'(X)=\{ \alpha_i \mid 1\leq i\leq s \}$.
For any $z\in X$, $X_i(z)$ is a $t_i:=t+1-s^*=(s+1)$-design in $S^{d-2}$. 
Then the degrees $s_{i,j}=|A(X_i(z),X_j(z))|$ satisfy $s_{i,j}\leq s$, and the following holds:
\begin{eqnarray*}
2s-2 \leq s+1&\Leftrightarrow & s \leq 3\\
&\Leftrightarrow & t = 2,4,6.
\end{eqnarray*}
In particular, if $t=4$, then $s_{i,j}+s_{j,k}-2\leq t_j$ holds,
i.e., the assumption (1) of Theorem~\ref{coherent} holds for all $i,j,k$.
By Corollary~\ref{triple}, we obtain the following result.
\begin{corollary}
Every tight $4$-design carries a triply regular association scheme. 
\end{corollary}
The same argument shows that a spherical $3$-design with degree $2$ i.e.,
a strongly regular graph with $a_1^*=0$ 
carries a triply regular association scheme. 
This is already known (see \cite{CGS}).

Let $X \subset S^{d-1}$ be a tight $(2s-1)$-design, and let 
$A'(X)=\{ \alpha_i \mid 1\leq i\leq s \}$ where $\alpha_s=-1$.
For any $z\in X$ and $i\neq s$, $X_i(z)$ is a $t_i:=t+1-s^*=(s+1)$-design in $S^{d-2}$. 

Then the degrees $s_{i,j}=|A(X_i(z),X_j(z))|$ satisfy $s_{i,j}\leq s-1$, and the following holds:
\begin{eqnarray*}
2s-4 \leq s+1&\Leftrightarrow & s \leq 5\\
                  &\Leftrightarrow & t = 1,3,5,7,9. 
\end{eqnarray*}
In particular, if $t=5,7$, then $s_{i,j}+s_{j,k}-2\leq t_j$ holds,
i.e., the assumption (1) of Theorem~\ref{coherent} holds for all $i,j,k$.
By Corollary~\ref{triple}, we obtain the following result. 
\begin{corollary}
Every tight $5$- or $7$-design carries a triply regular association scheme. 
\end{corollary}
The same argument shows that an antipodal spherical $3$-design with degree $3$ carries a triply regular association scheme i.e.,
subconstituents of a Taylor graph are strongly regular graphs.
This is already known (see \cite[Theorem 1.5.3]{BCN}).
\section{Derived designs of $Q$-polynomial association schemes}\label{sec:4}
The reader is referred to \cite{BI} for the basic information on $Q$-polynomial association schemes. 
The following lemma is used to prove Lemma~\ref{design}.
\begin{lemma}\label{as}
Let $\mathfrak{X}=(X,\{R_i\}_{i=0}^d)$ be a 
symmetric association scheme of class $d$.
Let $B_i=(p_{i,j}^k)$ be its $i$-th intersection matrix, 
and $Q=(q_j(i))$ be the second eigenmatrix of 
$\mathfrak{X}$.
Then 
\[
(Q^tB_i)(h,i)=\frac{k_iq_h(i)^2}{m_h}
\quad(0\leq h,i\leq d).
\]
\end{lemma}
\begin{proof}
See \cite[p.73 (4.2) and Theorem 3.5(i)]{BI}.
\end{proof}

The following lemma gives a property of derived designs of the embedding of a $Q$-polynomial association scheme into the first eigenspace.
\begin{lemma}\label{design}
Let $(X,\{R_i\}_{i=0}^s)$ be a $Q$-polynomial association scheme, 
and we identify $X$ as the image of the embedding into 
the first eigenspace by $E_1=\frac{1}{|X|}\sum\nolimits_{j=0}^s\theta_j^*A_j$.
Then, for $i\in\{1,\dots,s\}$ with $\theta_i^*\neq-\theta_0^*$,
the derived design $X_i(z)$ is a $2$-design in $S^{\theta_0^*-2}$
for any $z\in X$ if and only if $a_1^*(\ts_i+1)=0$.
\end{lemma} 
\begin{proof}
The angle set of $X_i(z)$ consists of
\[
\frac{\frac{\ts_k}{\ts_0}-\frac{{\ts_i}^2}{{\ts_0}^2}}{1
-(\frac{\ts_i}{\ts_0})^2}
=\frac{\ts_0\ts_k-{\ts_i}^2}{{\ts_0}^2-{\ts_i}^2}\quad
(0\leq k\leq s,\;p_{i,i}^k\neq0).
\]
Thus, Lemma~\ref{cha} implies that
$X_i(z)$ is a $2$-design in $S^{\theta_0^*-2}$ if and only if 
$$\sum\limits_{j=0}^sQ_k(\frac{\theta_0^*\theta_j^*-{\theta_i^*}^2}{{\theta_0^*}^2-{\theta_i^*}^2})p_{i,j}^i=0 \quad(k=1,2), $$
where $Q_k(x)$ is the Gegenbauer polynomial of degree $k$ in $S^{\theta_0^*-2}$.

Since $Q_1(x)=(\theta_0^*-1)x$, $\sum\limits_{j=0}^sp_{i,j}^j=k_i$ and
\begin{align}
\sum_{j=0}^s\ts_jp_{i,j}^i
&=(Q^tB_i)(1,i)
\nnexteq
\frac{k_iq_1(i)^2}{m_1}
&&\text{(by Lemma~\ref{as})}
\nnexteq
\frac{k_i{\ts_i}^2}{\ts_0},
\label{eq:1}
\end{align}
we have
\begin{align*}
\sum_{j=0}^s
Q_1(\frac{\theta_0^*\theta_j^*-{\theta_i^*}^2}{{\theta_0^*}^2-{\theta_i^*}^2})p_{i,j}^i
&=
\frac{\theta_0^*-1}{{\theta_0^*}^2-{\theta_i^*}^2}
\left(\ts_0
\sum_{j=0}^s\theta_j^* p_{i,j}^i
-{\theta_i^*}^2\sum_{j=0}^s p_{i,j}^i
\right)\\
&=0.
\end{align*}

Since $Q_2(x)=(\theta_0^*-1)x^2-1$, $\sum\limits_{j=0}^sp_{i,j}^j=k_i$, (\ref{eq:1}) and
\begin{align*}
\sum_{j=0}^s{\ts_j}^2 p_{i,j}^i
&=\sum_{j=0}^s (c_2^* q_2(i)+a_1^*q_1(i)
+b_0^*q_0(i))p_{i,j}^i
\nnexteq
c_2^*(Q^tB_i)(2,i)
+a_1^* \frac{k_i{\ts_i}^2}{\ts_0}
+\ts_0 k_i
&&\text{(by (\ref{eq:1}))}
\nnexteq
c_2^*\frac{k_iq_2(i)^2}{m_2}
+k_i(\frac{a_1^*{\ts_i}^2}{\ts_0}+\ts_0)
&&\text{(by Lemma~\ref{as})}
\nnexteq
k_i\left(
\frac{((\ts_i-a_1^*)\ts_i-\ts_0)^2}{
(\ts_0-a_1^*)\ts_0-\ts_0}
+\frac{a_1^*{\ts_i}^2}{\ts_0}+\ts_0\right),
\end{align*}
we have
\begin{align*}
\sum_{j=0}^s Q_2(\frac{\ts_0\ts_j-{\ts_i}^2}{
{\ts_0}^2-{\ts_i}^2})p_{i,j}^i
&=
\frac{\theta_0^*-1}{({\theta_0^*}^2-{\theta_i^*}^2)^2}
({\ts_0}^2\sum_{j=0}^s {\ts_j}^2p_{i,j}^i
-2\ts_0{\ts_i}^2\sum_{j=0}^s \ts_j p_{i,j}^i
+{\ts_i}^4\sum_{j=0}^s p_{i,j}^i)
-k_i
\nexteq
\frac{k_ia_1^*(\ts_i+1)^2\ts_0}{(\ts_0+\ts_i)^2(\ts_0-a_1^*-1)}.
\end{align*}
Therefore $X_i(z)$ is a $2$-design in $S^{\theta_0^*-2}$ if and only if $a_1^*(\ts_i+1)=0$.
\end{proof}
\section{Real mutually unbiased bases}\label{sec:5}
\begin{definition}
Let $M=\{M_i\}_{i=1}^f$ be a collection of orthonormal bases of $\mathbb{R}^d$.
$M$ is called  real mutually unbiased bases (MUB) if any two vectors $x$ and $y$ from different bases
satisfy $\langle x,y \rangle =\pm 1/\sqrt{d}$.
\end{definition}

It is known that the number $f$ of real mutually unbiased bases in $\mathbb{R}^d$ can be at most $d/2+1$.  
We call $M$ a maximal MUB if this upper bound is attained. 
Constructions of maximal MUB are known only for $d=2^{m+1}$, $m$ odd \cite{CCKS}.
Throughout this section, we assume $M=\{M_i\}_{i=1}^f$ is an MUB,  put $X^{(i)}=M_i \cup (-M_i)$ and $X=M \cup (-M)$. 
The angle set of $X$ is 
$$ A'(X)=\{\frac{1}{\sqrt{d}},0,-\frac{1}{\sqrt{d}},-1 \} .$$
We set 
$$\alpha_0=1,\quad\alpha_1=\frac{1}{\sqrt{d}},\quad\alpha_2=0,\quad \alpha_3=-\frac{1}{\sqrt{d}},\quad\alpha_4=-1 ,$$
and we define $R_k=\{(x,y)\in X\times X \mid \langle x,y\rangle=\alpha_k\}$.

Since $X^{(i)}$ is a spherical $3$-design in $S^{d-1}$
for any $i\in\{1,\ldots,f\}$,
$X$ is also a spherical $3$-design in $S^{d-1}$.
It is shown in \cite{LMO} that 
 $(X,\{R_k\}_{k=0}^4)$ is a $Q$-polynomial association scheme with $a_1^*=0$. 
$X$ is imprimitive and the set $\{X^{(1)},\ldots,X^{(f)}\}$ is a system of imprimitivity with respect to the equivalence relation $R_0\cup R_2\cup R_4$. 

By Lemma~\ref{design}, for any $z\in X$ the derived design $X_i=X_i(z)$ is a $t_i=2$-design in $S^{d-2}$.
We define $s_{i,j}=|A(X_i,X_j)|$. 
Then the matrix $(s_{i,j})_{\substack{1\leq i \leq 3\\1\leq j \leq 3}}$ is
$$\left(\begin{array}{ccc}
         3 & 2 & 3\\
         2 & 1 & 2\\
         3 & 2 & 3
\end{array}\right) .$$

If $s_{i,j}+s_{j,k}-2\leq 2$, that is, when
\begin{align*}
(i,j,k)\in&\{
(1,2,1), (1,2,2), (1,2,3), (2,1,2), (2,2,1), (2,2,2),\\
&(2,2,3), (2,3,2), (3,2,1), (3,2,2), (3,2,3)\},
\end{align*}
then the assumption (1) of Theorem~\ref{coherent} holds.
We remark that $X_2$ is in fact a $3$-design because $X_2$ is a cross polytope in $\mathbb{R}^{d-1}$,  
but this fact does not improve the proof.

The following Lemma is used to determine intersection numbers of derived designs obtained from MUB.
\begin{lemma}\label{mubsub}
We define $X_i(x,\alpha)=\{w\in X_i \mid \langle x,w\rangle=\alpha\}$, and $X_i(x,\alpha;y,\beta)=X_i(x,\alpha)\cap X_i(y,\beta)$.
Then the following equalities hold:
\begin{enumerate}
\item $X_i(x,-\alpha)=X_i(-x,\alpha)$,
\item $-X_i(x,\alpha)=X_{4-i}(x,-\alpha)$,
\item $|X_i(x,\alpha;y,\beta)|=|X_i(-x,-\alpha;y,\beta)|=|X_i(x,\alpha;-y,-\beta)|=|X_{4-i}(x,-\alpha;y,-\beta)|$.
\end{enumerate}
\end{lemma}
\begin{proof}
(1) and (2) are immediate from the definition.

By (1), 
$X_i(x,\alpha;y,\beta)=X_i(-x,-\alpha;y,\beta)=
X_i(x,\alpha;-y,-\beta)$ holds.
By (2), $-X_i(x,\alpha;y,\beta)=X_{4-i}(x,-\alpha;y,-\beta)$ holds.  
This proves (3).
\end{proof}

If $s_{i,j}+s_{j,k}-3=2$, that is, when
\begin{align}
(i,j,k)\in\{
(1,1,2), (1,3,2),(2,1,1),(2,1,3), 
      (2,3,1), (2,3,3),(3,1,2),(3,3,2)\},\label{MUB2s-3}
\end{align}
Lemma~\ref{mubsub} implies that the intersection numbers on $X_j(z)$ for $x\in X_i(z)$, $y\in X_k(z)$ are determined by the intersection numbers on $X_1(z)$ for $x'\in X_1(z)$, $y'\in X_2(z)$.
And the intersection numbers $p_{\alpha_{1,1}^2,\alpha_{1,2}^1}^1(x,y)$, $p_{\alpha_{1,1}^2,\alpha_{1,2}^3}^1(x,y)$ for $x,y\in X_1(z)$ are uniquely determined by $\gamma=\langle x,y\rangle$ as follows:
\begin{align*}
p_{\alpha_{1,1}^2,\alpha_{1,2}^1}^1(x,y)=\begin{cases}
     \frac{d}{2}-1 & \text{ if }\ \langle x,y\rangle=\alpha_{1,2}^1, \\
     \frac{d}{2}   & \text{ if }\ \langle x,y\rangle=\alpha_{1,2}^3, 
\end{cases} \quad p_{\alpha_{1,1}^2,\alpha_{1,2}^3}^1(x,y)=\begin{cases}
     \frac{d}{2} & \text{ if }\ \langle x,y\rangle=\alpha_{1,2}^1, \\
     \frac{d}{2}-1   & \text{ if }\ \langle x,y\rangle=\alpha_{1,2}^3. 
\end{cases}
\end{align*}
These numbers are independent of $z\in X$.
Hence the assumption (2) of Theorem~\ref{coherent} holds for $(i,j,k)$ in (\ref{MUB2s-3}).

If $s_{i,j}+s_{j,k}-4= 2$, that is, when
\begin{align}
(i,j,k)\in\{(1,1,1), (1,1,3),(1,3,1),(1,3,3),(3,1,1), (3,1,3),(3,3,1),(3,3,3)\},\label{MUB2s-4}
\end{align}
Lemma~\ref{mubsub} implies that the intersection numbers on $X_j(z)$ for $x\in X_i(z)$, $y\in X_k(z)$ are determined by the intersection numbers on $X_1(z)$ for $x'\in X_1(z)$, $y'\in X_1(z)$.
And the intersection numbers $\{p_{\alpha,\beta}^1(x,y) \mid \alpha=\alpha_{1,1}^2\ \text{or}\ \beta=\alpha_{1,1}^2\} $ are given in Table~\ref{tb:t1}.
These numbers are independent of $z\in X$. 
Hence the assumption (3) of Theorem~\ref{coherent} holds for $(i,j,k)$ in (\ref{MUB2s-4}).   
By Corollary~\ref{triple}, we obtain the following result.
\begin{corollary}
Every MUB carries a triply regular association scheme. 
\end{corollary}

\begin{table}[h]
\begin{center}
\caption{the values of $p_{\alpha,\beta}^1(x,y)$, where $x\in X_1$, $y\in X_1$}
\begin{tabular}{|c|c|} \hline 
$(\alpha,\beta)$ & $p_{\alpha,\beta}^1(x,y)$ \\ \hline
$(\alpha_{1,1}^2,\alpha_{1,1}^2)$ & 
$\left\{\begin{array}{ll}
     0 & \mbox{if}\ \langle x,y\rangle=\alpha_{1,1}^1 \\
     d-2 & \mbox{if}\ \langle x,y\rangle=\alpha_{1,1}^2 \\
     0 & \mbox{if}\ \langle x,y\rangle=\alpha_{1,1}^3 
\end{array}\right.$ \\ \hline
$\begin{array}{c}(\alpha_{1,1}^2,\alpha_{1,1}^1),\\(\alpha_{1,1}^1,\alpha_{1,1}^2)\end{array}$ & $\left\{\begin{array}{ll}
     \frac{d+\sqrt{d}}{2}-1 & \mbox{if}\ \langle x,y\rangle=\alpha_{1,1}^1 \\
     0 & \mbox{if}\ \langle x,y\rangle=\alpha_{1,1}^2 \\
     \frac{d+\sqrt{d}}{2} & \mbox{if}\ \langle x,y\rangle=\alpha_{1,1}^3 
\end{array}\right.$ \\ \hline
$\begin{array}{c}(\alpha_{1,1}^2,\alpha_{1,1}^3),\\(\alpha_{1,1}^3,\alpha_{1,1}^2)\end{array}$ & $\left\{\begin{array}{ll}
     \frac{d-\sqrt{d}}{2} & \mbox{if}\ \langle x,y\rangle=\alpha_{1,1}^1 \\
     0 & \mbox{if}\ \langle x,y\rangle=\alpha_{1,1}^2 \\
     \frac{d-\sqrt{d}}{2}-1 & \mbox{if}\ \langle x,y\rangle=\alpha_{1,1}^3 
\end{array}\right.$\\ \hline
\end{tabular}\label{tb:t1}
\end{center}
\end{table}
\section{Linked systems of symmetric designs}\label{sec:6}
\begin{definition}
Let $(\Omega_i,\Omega_j,I_{i,j})$ be an incidence structure satisfying 
$\Omega_i\cap \Omega_j=\emptyset$, $I_{j,i}^t=I_{i,j}$ for any distinct integers $i,j \in \{1,\dots,f\}$. 
We put $\Omega=\bigcup_{i=1}^f\Omega_i$,
$I=\bigcup_{i\neq j} I_{i,j}$.
$(\Omega, I)$ is called a linked system of symmetric $(v,k,\lambda)$ designs if the following conditions hold:
\begin{enumerate}
\item for any distinct integers $i,j\in\{1,\dots,f\}$, $(\Omega_i,\Omega_j,I_{i,j})$ is a symmetric $(v,k,\lambda)$ design,
\item for any distinct integers $i,j,l \in\{1,\dots,f\}$, and for any $x\in \Omega_i, y\in \Omega_j$, the number of $z\in \Omega_l$ incident with 
both $x$ and $y$ depends only on whether $x$ and $y$ are incident or not, and does not depend on $i,j,l$.
\end{enumerate}
\end{definition}
We define the integers $\sigma,\tau$ by
\[
|\{z\in \Omega_l \mid (x,z)\in I_{i,l}, (y,z)\in I_{j,l}\}|=
\begin{cases}\sigma&\text{ if }\ (x,y)\in I_{i,j},\\
\tau&\text{ if }\ (x,y)\not\in I_{i,j},
\end{cases}
\]
where $i,j,l \in\{1,\dots,f\}$ are distinct and $x\in \Omega_i$, $y\in \Omega_j$.

By \cite[Theorem 1]{C}, we may assume that 
$$ \sigma=\frac{1}{v}(k^2-\sqrt{n}(v-k)),\quad \tau=\frac{k}{v}(k+\sqrt{n}),  $$
where $n=k-\lambda$.
It is easy to see that $(\Omega,\{R_i\}_{i=0}^3)$
is a $3$-class association scheme, where
\begin{align*}
R_0&= \{(x,x) \mid x\in \Omega\} ,\\
R_1&= \{(x,y)\mid x\in \Omega_i,y\in \Omega_j, (x,y)\in I_{i,j}\ \mbox{for some}\ i\neq j\} ,\\
R_2&= \{(x,y)\mid x,y\in \Omega_i, x\neq y \mbox{ for some}\ i\} ,\\
R_3&= \{(x,y)\mid x\in \Omega_i,y\in \Omega_j, (x,y)\not\in I_{i,j}\ \mbox{for some}\ i\neq j\}.
\end{align*} 

We note that the second eigenmatrix $Q$ is given in \cite{M}
as follows:
\[
Q=\left(\begin{array}{cccc}
 1 & v-1 & (f-1)(v-1) & f-1 \\
 1 & -\sqrt{\frac{(v-1)(v-k)}{k}} & \sqrt{\frac{(v-1)(v-k)}{k}} & -1 \\
 1 & -1 & -f+1 & f-1 \\
 1 & \sqrt{\frac{(v-1)k}{v-k}} & -\sqrt{\frac{(v-1)k}{v-k}} & -1
 \end{array}\right) ,
\]
and hence the Krein matrix
$B_1^*=(q_{1,j}^k)_{\substack{0\leq j\leq 3\\0\leq k\leq 3}}$
is given as follows:
\[
B_1^*=\left(\begin{array}{cccc}
 0 & 1 & 0 & 0 \\
 v-1 & \frac{k(v-k)(v-2)+(f-1)(2k-v)\sqrt{k(v-k)(v-1)}}{fk(v-k)} & \frac{k(v-k)(v-2)+(v-2k)\sqrt{k(v-k)(v-1)}}{fk(v-k)} & 0 \\
 0 & \frac{(f-1)(k(v-k)(v-2)+(v-2k)\sqrt{k(v-k)(v-1)})}{fk(v-k)} & \frac{(f-1)k(v-k)(v-2)+(2k-v)\sqrt{k(v-k)(v-1)}}{fk(v-k)} & v-1 \\
 0 & 0 & 1 & 0
 \end{array}\right).
\]
Therefore $(\Omega,\{R_i\}_{i=0}^3)$ is a $Q$-polynomial association scheme. 
$(\Omega,\{R_i\}_{i=0}^3)$ is imprimitive and the set $\{\Omega_1,\ldots,\Omega_f\}$ is a system of imprimitivity with respect to the equivalence relation $R_0\cup R_2$. 

In the rest of this section, we assume that $a_1^*=0$ i.e., 
$f=1+\frac{(v-2)\sqrt{k(v-k)}}{(v-2k)\sqrt{v-1}}$.
Examples of linked symmetric designs satisfying this
assumption are known for  $(v,k,\lambda)=(2^{2m},2^{2m-1}-2^{m-1},2^{2m-2}-2^{m-1})$ with $f=2^{2m-1}$ for any $m>1$ \cite{C}. 
 
Let $X$ be the embedding of $\Omega$ into the first eigenspace. 
The angle set of $X$ is 
$$A'(X)=\{\frac{\theta_k^*}{\theta_0^*}\mid 1\leq k\leq3\},$$
and we set 
$\alpha_k=\theta_k^*/\theta_0^*$.
We consider the derived design $X_i(z)$ for $z\in X$.
By $a_1^*=0$, Lemma~\ref{design} implies $X_i(z)$ is a $2$-design in $S^{v-3}$.
We define $s_{i,j}=|A'(X_i(z),X_j(z))|$. 
Then the matrix $(s_{i,j})_{\substack{1\leq i \leq 3\\1\leq j \leq 3}}$ is
$$\left(\begin{array}{ccc}
         3 & 2 & 3\\
         2 & 1 & 2\\
         3 & 2 & 3
\end{array}\right) .$$
Since $\{\Omega_1,\ldots,\Omega_f\}$ is a system of imprimitivity, 
we obtain Table~\ref{tb:t2}, Table~\ref{tb:t3}.

If $s_{i,j}+s_{j,l}-2\leq 2$, that is, when
\begin{align*}
(i,j,l)\in&\{
(1,2,1), (1,2,2), (1,2,3), (2,1,2), (2,2,1), (2,2,2),\\
&(2,2,3), (2,3,2), (3,2,1), (3,2,2), (3,2,3)\},
\end{align*}
then the assumption (1) of Theorem~\ref{coherent} holds.

If $s_{i,j}+s_{j,l}-3= 2$, that is, when
\begin{align}
(i,j,l)\in\{(1,1,2), (1,3,2), (2,1,1), (2,1,3), (2,3,1), (2,3,3), (3,1,2), (3,3,2)\}, \label{LSD2s-3}
\end{align}
Table~\ref{tb:t2} implies that the numbers $p_{\alpha_{i,j}^2,\alpha_{j,l}^1}^j(x,y)$ or $p_{\alpha_{i,j}^1,\alpha_{j,l}^2}^j(x,y)$ are independent of
$z\in X$ and $(x,y)\in X_i(z)\times X_l(z)$ with
$\gamma=\langle x,y\rangle$.
Hence the assumption (2) of Theorem~\ref{coherent} holds for $(i,j,l)$ in (\ref{LSD2s-3}). 

If $s_{i,j}+s_{j,l}-4= 2$, that is, when
\begin{align}
(i,j,l)\in\{(1,1,1), (1,1,3), (1,3,1), (1,3,3), (3,1,1), (3,1,3), (3,3,1), (3,3,3)\}, \label{LSD2s-4}
\end{align}
Table~\ref{tb:t3} implies the numbers $p_{\alpha_{i,j}^2,\alpha_{j,l}^2}^j(x,y)$, $p_{\alpha_{i,j}^2,\alpha_{j,l}^1}^j(x,y)$ and $p_{\alpha_{i,j}^1,\alpha_{j,l}^2}^j(x,y)$ are independent of
$z\in X$ and $(x,y)\in X_i(z)\times X_l(z)$ with
$\gamma=\langle x,y\rangle$.
Hence the assumption (3) of Theorem~\ref{coherent} holds for $(i,j,l)$ in (\ref{LSD2s-4}).
By Corollary~\ref{triple}, we obtain the following result.
\begin{corollary}
Every linked system of symmetric design satisfying $f=1+\frac{(v-2)\sqrt{k(v-k)}}{(v-2k)\sqrt{v-1}}$ carries a triply regular association scheme. 
\end{corollary}

\begin{table}[h]
\begin{center}
\caption{the values of $p_{\alpha,\beta}^j(x,y)$, where $x\in X_i(z)$, $y\in X_l(z)$}
\begin{tabular}{|@{}c@{}|c|c||@{}c@{}|c|c|}\hline
$(i,j,l)$&$(\alpha,\beta)$&$p_{\alpha,\beta}^j(x,y)$ &$(i,j,l)$&$(\alpha,\beta)$&$p_{\alpha,\beta}^j(x,y)$\\ \hline
$(1,1,2)$ & $(\alpha_{1,1}^2,\alpha_{1,2}^1)$ & 
$\left\{\begin{array}{ll}
     \lambda-1 &  \langle x,y\rangle=\alpha_{1,2}^1 \\
     \lambda & \langle x,y\rangle=\alpha_{1,2}^3
\end{array}\right.$ &$(2,1,1)$ & $(\alpha_{2,1}^1,\alpha_{1,1}^2)$ & $\left\{\begin{array}{ll}
     \lambda-1 & \langle x,y\rangle=\alpha_{2,1}^1 \\
     \lambda & \langle x,y\rangle=\alpha_{2,1}^3
\end{array}\right.$ \\ \hline
$(1,3,2)$ & $(\alpha_{1,3}^2,\alpha_{3,2}^1)$ & 
$\left\{\begin{array}{ll}
     k-\lambda &  \langle x,y\rangle=\alpha_{1,2}^1 \\
     k-\lambda & \langle x,y\rangle=\alpha_{1,2}^3
\end{array}\right.$ &$(2,3,1)$ & $(\alpha_{2,3}^1,\alpha_{3,1}^2)$ & $\left\{\begin{array}{ll}
     k-\lambda & \langle x,y\rangle=\alpha_{2,1}^1 \\
     k-\lambda & \langle x,y\rangle=\alpha_{2,1}^3
\end{array}\right.$\\ \hline
$(3,1,2)$ & $(\alpha_{3,1}^2,\alpha_{1,2}^1)$ & 
$\left\{\begin{array}{ll}
     \lambda &  \langle x,y\rangle=\alpha_{3,2}^1 \\
     \lambda & \langle x,y\rangle=\alpha_{3,2}^3
\end{array}\right.$ &$(2,1,3)$ & $(\alpha_{2,1}^1,\alpha_{1,3}^2)$ & $\left\{\begin{array}{ll}
     \lambda & \langle x,y\rangle=\alpha_{2,3}^1 \\
     \lambda & \langle x,y\rangle=\alpha_{2,3}^3
\end{array}\right.$\\ \hline
$(3,3,2)$ & $(\alpha_{3,3}^2,\alpha_{3,2}^1)$ & 
$\left\{\begin{array}{ll}
     k-\lambda-1 &  \langle x,y\rangle=\alpha_{3,2}^1 \\
     k-\lambda & \langle x,y\rangle=\alpha_{3,2}^3
\end{array}\right.$ &$(2,3,3)$ & $(\alpha_{2,3}^1,\alpha_{3,3}^2)$ & $\left\{\begin{array}{ll}
     k-\lambda-1 & \langle x,y\rangle=\alpha_{2,3}^1 \\
     k-\lambda & \langle x,y\rangle=\alpha_{2,3}^3
\end{array}\right.$\\ \hline
\end{tabular}\label{tb:t2}
\end{center}
\end{table}

\begin{table}[t]
\begin{center}
\caption{the values of $p_{\alpha,\beta}^j(x,y)$, where $x\in X_i(z)$, $y\in X_l(z)$}
\begin{tabular}{|@{}c@{}|c|c||@{}c@{}|c|c|}\hline
$(i,j,l)$&$(\alpha,\beta)$&$p_{\alpha,\beta}^j(x,y)$ &$(i,j,l)$&$(\alpha,\beta)$&$p_{\alpha,\beta}^j(x,y)$\\ \hline
 & $(\alpha_{1,1}^2,\alpha_{1,1}^2)$ & 
$\left\{\begin{array}{ll}
     0 & \langle x,y\rangle=\alpha_{1,1}^1 \\
     k-2 & \langle x,y\rangle=\alpha_{1,1}^2 \\
     0 & \langle x,y\rangle=\alpha_{1,1}^3
\end{array}\right.$ & &$(\alpha_{1,3}^2,\alpha_{3,3}^2)$ &$\left\{\begin{array}{ll}
     0 & \langle x,y\rangle=\alpha_{1,3}^1 \\
     v-k-1 & \langle x,y\rangle=\alpha_{1,3}^2 \\
     0 & \langle x,y\rangle=\alpha_{1,3}^3
\end{array}\right.$ \\ \cline{2-3} \cline{5-6}
$(1,1,1)$ & $(\alpha_{1,1}^2,\alpha_{1,1}^1)$ & 
$\left\{\begin{array}{ll}
     \sigma-1 &  \langle x,y\rangle=\alpha_{1,1}^1 \\
     0 & \langle x,y\rangle=\alpha_{1,1}^2 \\
     \sigma & \langle x,y\rangle=\alpha_{1,1}^3
\end{array}\right.$ &$(1,3,3)$ & $(\alpha_{1,3}^2,\alpha_{3,3}^1)$ & $\left\{\begin{array}{ll}
     k-\tau & \langle x,y\rangle=\alpha_{1,3}^1 \\
     0 & \langle x,y\rangle=\alpha_{1,3}^2 \\
     k-\tau & \langle x,y\rangle=\alpha_{1,3}^3
\end{array}\right.$ \\ \cline{2-3} \cline{5-6}
 & $(\alpha_{1,1}^1,\alpha_{1,1}^2)$ & 
$\left\{\begin{array}{ll}
     \sigma-1 & \langle x,y\rangle=\alpha_{1,1}^1 \\
     0 & \langle x,y\rangle=\alpha_{1,1}^2 \\
     \sigma & \langle x,y\rangle=\alpha_{1,1}^3
\end{array}\right.$ & &$(\alpha_{1,3}^1,\alpha_{1,3}^2)$ &$\left\{\begin{array}{ll}
     k-\sigma-1 & \langle x,y\rangle=\alpha_{1,3}^1 \\
     0 & \langle x,y\rangle=\alpha_{1,3}^2 \\
     k-\sigma & \langle x,y\rangle=\alpha_{1,3}^3
\end{array}\right.$  \\ \hline
 & $(\alpha_{1,1}^2,\alpha_{1,3}^2)$ & 
$\left\{\begin{array}{ll}
     0 & \langle x,y\rangle=\alpha_{1,3}^1 \\
     k-1 & \langle x,y\rangle=\alpha_{1,3}^2 \\
     0 & \langle x,y\rangle=\alpha_{1,3}^3
\end{array}\right.$ & &$(\alpha_{3,1}^2,\alpha_{1,3}^2)$ &$\left\{\begin{array}{ll}
     0 & \langle x,y\rangle=\alpha_{3,3}^1 \\
     k & \langle x,y\rangle=\alpha_{3,3}^2 \\
     0 & \langle x,y\rangle=\alpha_{3,3}^3
\end{array}\right.$ \\ \cline{2-3} \cline{5-6}
$(1,1,3)$ & $(\alpha_{1,1}^2,\alpha_{1,3}^1)$ & 
$\left\{\begin{array}{ll}
     \tau-1 &  \langle x,y\rangle=\alpha_{1,3}^1 \\
     0 & \langle x,y\rangle=\alpha_{1,3}^2 \\
     \tau & \langle x,y\rangle=\alpha_{1,3}^3
\end{array}\right.$ &$(3,1,3)$ & $(\alpha_{3,1}^2,\alpha_{1,3}^1)$ & $\left\{\begin{array}{ll}
     \tau & \langle x,y\rangle=\alpha_{3,3}^1 \\
     0 & \langle x,y\rangle=\alpha_{3,3}^2 \\
     \tau & \langle x,y\rangle=\alpha_{3,3}^3
\end{array}\right.$ \\ \cline{2-3} \cline{5-6}
 & $(\alpha_{1,1}^1,\alpha_{1,3}^2)$ & 
$\left\{\begin{array}{ll}
     \sigma & \langle x,y\rangle=\alpha_{1,3}^1 \\
     0 & \langle x,y\rangle=\alpha_{1,3}^2 \\
     \sigma & \langle x,y\rangle=\alpha_{1,3}^3
\end{array}\right.$ & &$(\alpha_{3,1}^1,\alpha_{1,3}^2)$ &$\left\{\begin{array}{ll}
     \tau & \langle x,y\rangle=\alpha_{3,3}^1 \\
     0 & \langle x,y\rangle=\alpha_{3,3}^2 \\
     \tau & \langle x,y\rangle=\alpha_{3,3}^3
\end{array}\right.$  \\ \hline
 & $(\alpha_{1,3}^2,\alpha_{3,1}^2)$ & 
$\left\{\begin{array}{ll}
     0 & \langle x,y\rangle=\alpha_{1,1}^1 \\
     v-k & \langle x,y\rangle=\alpha_{1,1}^2 \\
     0 & \langle x,y\rangle=\alpha_{1,1}^3
\end{array}\right.$ & &$(\alpha_{3,3}^2,\alpha_{3,1}^2)$ &$\left\{\begin{array}{ll}
     0 & \langle x,y\rangle=\alpha_{3,1}^1 \\
     v-k-1 & \langle x,y\rangle=\alpha_{3,1}^2 \\
     0 & \langle x,y\rangle=\alpha_{3,1}^3
\end{array}\right.$ \\ \cline{2-3} \cline{5-6}
$(1,3,1)$ & $(\alpha_{1,3}^2,\alpha_{3,1}^1)$ & 
$\left\{\begin{array}{ll}
     k-\sigma &  \langle x,y\rangle=\alpha_{1,1}^1 \\
     0 & \langle x,y\rangle=\alpha_{1,1}^2 \\
     k-\sigma & \langle x,y\rangle=\alpha_{1,1}^3
\end{array}\right.$ &$(3,3,1)$ & $(\alpha_{3,3}^2,\alpha_{3,1}^1)$ & $\left\{\begin{array}{ll}
     k-\tau-1 & \langle x,y\rangle=\alpha_{3,1}^1 \\
     0 & \langle x,y\rangle=\alpha_{3,1}^2 \\
     k-\tau & \langle x,y\rangle=\alpha_{3,1}^3
\end{array}\right.$ \\ \cline{2-3} \cline{5-6}
 & $(\alpha_{1,3}^1,\alpha_{3,1}^2)$ & 
$\left\{\begin{array}{ll}
     k-\sigma & \langle x,y\rangle=\alpha_{1,1}^1 \\
     0 & \langle x,y\rangle=\alpha_{1,1}^2 \\
     k-\sigma & \langle x,y\rangle=\alpha_{1,1}^3
\end{array}\right.$ & &$(\alpha_{3,3}^1,\alpha_{3,1}^2)$ &$\left\{\begin{array}{ll}
     k-\tau & \langle x,y\rangle=\alpha_{3,1}^1 \\
     0 & \langle x,y\rangle=\alpha_{3,1}^2 \\
     k-\tau & \langle x,y\rangle=\alpha_{3,1}^3
\end{array}\right.$  \\ \hline
 & $(\alpha_{3,1}^2,\alpha_{1,1}^2)$ & 
$\left\{\begin{array}{ll}
     0 & \langle x,y\rangle=\alpha_{3,1}^1 \\
     k-1 & \langle x,y\rangle=\alpha_{3,1}^2 \\
     0 & \langle x,y\rangle=\alpha_{3,1}^3
\end{array}\right.$ & &$(\alpha_{3,3}^2,\alpha_{3,3}^2)$ &$\left\{\begin{array}{ll}
     0 & \langle x,y\rangle=\alpha_{3,3}^1 \\
     v-k-2 & \langle x,y\rangle=\alpha_{3,3}^2 \\
     0 & \langle x,y\rangle=\alpha_{3,3}^3
\end{array}\right.$ \\ \cline{2-3} \cline{5-6}
$(3,1,1)$ & $(\alpha_{3,1}^2,\alpha_{1,1}^1)$ & 
$\left\{\begin{array}{ll}
     \sigma &  \langle x,y\rangle=\alpha_{3,1}^1 \\
     0 & \langle x,y\rangle=\alpha_{3,1}^2 \\
     \sigma & \langle x,y\rangle=\alpha_{3,1}^3
\end{array}\right.$ &$(3,3,3)$ & $(\alpha_{3,3}^2,\alpha_{3,3}^1)$ & $\left\{\begin{array}{ll}
     k-\tau-1 & \langle x,y\rangle=\alpha_{3,3}^1 \\
     0 & \langle x,y\rangle=\alpha_{3,3}^2 \\
     k-\tau & \langle x,y\rangle=\alpha_{3,3}^3
\end{array}\right.$ \\ \cline{2-3} \cline{5-6}
 & $(\alpha_{3,1}^1,\alpha_{1,1}^2)$ & 
$\left\{\begin{array}{ll}
     \tau-1 & \langle x,y\rangle=\alpha_{3,1}^1 \\
     0 & \langle x,y\rangle=\alpha_{3,1}^2 \\
     \tau & \langle x,y\rangle=\alpha_{3,1}^3
\end{array}\right.$ & &$(\alpha_{3,3}^1,\alpha_{3,3}^2)$ &$\left\{\begin{array}{ll}
     k-\tau-1 & \langle x,y\rangle=\alpha_{3,3}^1 \\
     0 & \langle x,y\rangle=\alpha_{3,3}^2 \\
     k-\tau & \langle x,y\rangle=\alpha_{3,3}^3
\end{array}\right.$  \\ \hline
\end{tabular}\label{tb:t3}
\end{center}
\end{table}
\section*{Acknowledgements}
The author would like to thank Professor Akihiro Munemasa 
for helpful discussions.
This work was supported by Grant-in-Aid for JSPS Fellows.

\end{document}